\documentclass[notitlepage, twocolumn, letterpaper, prx, longbibliography, 10pt]{revtex4-2}
\usepackage[utf8]{inputenc}
\usepackage{amsmath}
\usepackage{amssymb}
\usepackage{graphicx}
\usepackage{hyperref}
\usepackage{subcaption}
\usepackage[percent]{overpic}
\usepackage{pdfpages}
\usepackage[normalem]{ulem}
\usepackage[switch,columnwise]{lineno}
\usepackage[margin=1 in]{geometry}
\usepackage[object=vectorian]{pgfornament} % fancy hrules

\newcommand{\dd}[1]{\mathrm{d}#1}
\newcommand{\D}[1]{\Delta#1}

\makeatletter
\AtBeginDocument{\let\LS@rot\@undefined}
\makeatother

\graphicspath{{images/}}

\begin{document}

%************* Front matter ******************
\title{Equilibrium Moment Analysis of It\^o SDEs}
%\title{Modeling and analysis of systems with nonlinear functional dependence on random quantities}
%\title{Numerical simulation of nonlinear dependence on random quantities}

\author{David Sabin-Miller}
%\email[]{\mbox{davidsabinmiller@u.northwestern.edu}}
\email{dasami@umich.edu}
\affiliation{Center for the Study of Complex Systems, University of Michigan, Ann Arbor, MI, USA}
\author{Daniel M.~Abrams}
\email{dmabrams@northwestern.edu}
\affiliation{Department of Engineering Sciences and Applied Mathematics, Northwestern University, Evanston, IL, USA}
\affiliation{Department of Physics and Astronomy, Northwestern University, Evanston, IL, USA}
%\affiliation{Northwestern Institute for Complex Systems, Northwestern University, Evanston, IL, USA}

\begin{abstract}
    Stochastic differential equations have proved to be a valuable governing framework for many real-world systems which exhibit ``noise'' or randomness in their evolution.  One quality of interest in such systems is the shape of their equilibrium probability distribution, if such a thing exists. In some cases a straightforward integral equation may yield this steady-state distribution, but in other cases the equilibrium distribution exists and yet that integral equation diverges. Here we establish a new equilibrium-analysis technique based on the logic of finite-timestep simulation which allows us to glean information about the equilibrium regardless---in particular, a relationship between the raw moments of the equilibrium distribution. We utilize this technique to extract information about one such equilibrium resistant to direct definition.
\end{abstract}

\maketitle

%\section{Equilibrium Moment Analysis of It\^o SDEs}

% % DMA version of first 2 pars
% We now shift our focus to a technique of equilibrium analysis for It\^o SDEs. This is not directly tied to the work outlined in the previous section, but is of utility in analyzing a key feature of potential interest in any It\^o system, namely, its equilibrium distribution (if it exists).

% In some cases (including Eq.~\eqref{eq:nIto_form_of_cubic}) direct calculation of the steady-state distribution (as described in, e.g.,~\cite{gardiner2009stochastic}) requires integrals that fail to converge. We will now put forward a technique which avoids this problem and yields insight about this distribution. This technique may be applied to any It\^o system, even if the steady-state isn't directly calculable.
% %%%%%

% DSM version of first pars
%We now shift our focus from It\^o interpretation of generalized Langevin equations to a technique of equilibrium analysis for It\^o SDEs themselves. 
Equilibrium distributions are of considerable interest in any system where they exist. However, in some cases,  direct analytical calculation of the steady-state distribution (as described in, e.g.,~\cite{gardiner2009stochastic}) requires integrals that fail to converge. %Our interpretation of the ``stochastic cubic attractor'' from a forthcoming work \cite{DSM_nonlin_stochastic_interp} \eqref{eq:nIto_form_of_cubic} is of this type, and we will show that this technique yields insight into its structure.

% % for cases where an equilibrium distribution exists but is not analytically attainable. While not directly tied to our nonlinear-Langevin proposal, the technique we put forward is of utility in analyzing the nonlinear-Langevin cubic attractor example from earlier, which defies exact solution. 
% %%%%%

% % \dsm{I feel like there should be a clearer transition like this here, but it could use some smoothing into the next part}

% One key feature of potential interest in any It\^o system is its equilibrium distribution (if it exists). However, in some cases (including Eq.~\eqref{eq:nIto_form_of_cubic}) direct calculation of the steady-state distribution (as described in, e.g.,~\cite{gardiner2009stochastic}) requires integrals that fail to converge. We will now put forward a technique which avoids this problem and yields insight about this distribution. This technique may be applied to any It\^o system, even if the steady-state isn't directly calculable.

Suppose we seek to examine the equilibrium distribution (if it exists) of the autonomous It\^o SDE
\begin{align}
    \dd{x} = F(x) \dd{t} + G(x) \dd{W}.
\end{align}
We will use Euler-Maruyama numerical integration \cite{maruyama1955continuous} as a guide: in discrete time, we have
\begin{align}
    \D{x} = F(x) \D{t} + G(x) \eta \sqrt{\D{t}} , 
\end{align}
where $\eta \sim N(0,1)$. We may write the expression for the distribution of the new value $\xi = x + \D{x}$ from any previous position $x$:
\begin{align}
    \xi &\sim  N\left(x+F(x)\D{t}\ ,\ G(x)\sqrt{\D{t}}\right)\;, \nonumber \\
    P(\xi | x) &= \frac{1}{G(x)\sqrt{2 \pi \D{t}}}e^{\frac{-[\xi - x - F(x)\D{t}]^2}{2 G(x)^2 \D{t}}}\;.
\end{align}

%\dma{Probably not worth changing, but we did alrady use the symbol $P$ for something else on page 2 of the manuscript.  I can't think of a change that would improve things though, so let's leave it as is.}
%\dsm{I was using that to mean "probability of" rather than a specific function}
%\subsection{B. PDF of equilibrium}
Given this probability density function (PDF) for the outcome of a single step from any initial position $x$, we may write an expression for the evolution of the solution PDF from initial state $\rho_k(x)$ to subsequent state $\rho_{k+1}(x)$ a short time $\D{t}$ later:
\begin{equation*}
   \rho_{k+1}(\xi ) = \int\limits_{-\infty}^{\infty} P(\xi|x) \rho_k(x) \dd{x}\;. 
\end{equation*}
At equilibrium, this operation leaves the distribution $\rho_k=\rho_{k+1}=\rho^*$ unchanged, i.e.,
\begin{equation}
    \rho^*(\xi) = \int\limits_{-\infty}^{\infty} P(\xi|x) \rho^*(x) \dd{x}\;.
    \label{eq:rho*_implicit}
    % \rho^*(\xi) &= \frac{1}{n\sigma \sqrt{2\pi} \dt^{\frac{\alpha}{n}-\beta}} \cdot \nonumber\\
    %               & \ \ \ \ \  \int\limits_{-\infty}^{\infty} \left(\xi-x\right)^\frac{1-n}{n} e^{\frac{-\left[\left(\xi-x\right)^\frac{1}{n}+\D{t}^\frac{\alpha}{n}x \right]^2}{2\sigma^2\D{t}^{2\left(\frac{\alpha}{n}-\beta\right)}}} \rho^*(x) \dd{x} \label{eq:rho*_implicit}
\end{equation}
\subsection{Second Moment Method}
Rather than attempt to solve this implicit integral equation for $\rho^*$ directly, we instead examine the second (raw) moment of the distribution $\mu_2$ by multiplying both sides of Eq.~\eqref{eq:rho*_implicit} by $\xi^2$ and integrating over all $\xi$:
\begin{align*}
   %& \textrm{var}(\rho^*) = 
    &\mu_2 = \int\limits_{-\infty}^{\infty}\xi^2\rho^*(\xi) \dd{\xi} 
= \int\limits_{-\infty}^{\infty}\xi^2 \left[\int\limits_{-\infty}^{\infty} \rho^*(x) P(\xi | x) \dd{x} \right] \dd{\xi}  \\
  % &=  \int\limits_{-\infty}^{\infty}\xi^2 
  %     \int\limits_{-\infty}^{\infty} \rho^*(x)
  %       \frac{1}{G(x)\sqrt{2 \pi \D{t}}}e^{\frac{-[\xi - x - F(x)\D{t}]^2}{2 G(x)^2 \D{t}}} \dd{x} \ \dd{\xi}\\
   &= \int\limits_{-\infty}^{\infty}\rho^*(x) 
       \int\limits_{-\infty}^{\infty}\xi^2
         \frac{1}{G(x)\sqrt{2 \pi \D{t}}}e^{\frac{-[\xi - x - F(x)\D{t}]^2}{2 G(x)^2 \D{t}}} \dd{\xi}\ \dd{x}\;.
    %\label{eq:var_int}
\end{align*}
After swapping the order of integration\footnote{Changes in the order of integration will always be allowable for finite $\mu_2$.}, we observe that the inner integral over $\xi$ is of the form 
\begin{equation*}
    \frac{1}{s \sqrt{2 \pi}} \int\limits_{-\infty}^{\infty} u^2 e^{\frac{-(u-a)^2}{2s^2}} \dd{u} = a^2 + s^2
\end{equation*}
with $u=\xi$, $a = x+ F(x)\D{t}$, and $s = G(x) \sqrt{\D{t}}$. So we find
\begin{align*}
    \mu_2 %= &\int\limits_{-\infty}^{\infty}\rho^*(x) 
       %\left\{ \left[x+ F(x)\D{t}\right]^2 + G(x)^2\D{t}\right\} \dd{x} \\ 
    = &\int\limits_{-\infty}^{\infty} \rho^*(x) 
       \Big[ x^2+ 2xF(x)\D{t} + F(x)^2\D{t}^2 \\
       & \qquad \qquad \qquad \qquad \qquad\qquad + G(x)^2\D{t}\Big] \dd{x}\;.    
\end{align*}
Distributing the integral and subtracting $\mu_2$ from both sides (note that the integral of $x^2$ against $\rho^*$ is simply the definition of $\mu_2$), we find
\begin{align}
       0 &= \D{t} \int\limits_{-\infty}^{\infty}\rho^*(x) 
        \left[2xF(x) + G(x)^2\right] \dd{x}  \nonumber\\
          & \qquad+\D{t}^2\int\limits_{-\infty}^{\infty}\rho^*(x) 
        F(x)^2 \dd{x}\;, 
        \label{eq:consistency_Ito}
\end{align}
which should hold exactly for the equilibrium distribution(s) of any such It\^o system with small finite timestep $\Delta t$. We note, however, that getting to this point (swapping order of integration, and subtracting $\mu_2$ from both sides) includes the implicit assumption that $\mu_2$ is finite.

\subsection{Application to a Specific System}

We now look to apply this to a particular case, the ``cubic stochastic attractor'' from \cite{sabinmiller2024interpretation}:
\begin{align*}
    \frac{\dd{x}}{\dd{t}} = -X^3, \textrm{ where } X \sim N(x,\sigma)%-(x+\sigma \eta_t)^3
\end{align*}
which we argued is equivalent to the It\^o SDE 
\begin{align}
    \dd{x} &= (-x^3-3\sigma^2 x)\ \dd{t}  \nonumber \\
    & \qquad + \sqrt{15\sigma^6 + 36 \sigma^4 x^2 +9 \sigma^2 x^4}\ \dd{W}, \label{eq:nIto_form_of_cubic}
\end{align}
% So we argue that the system 
% \begin{align*}
%     \frac{\dd{x}}{\dd{t}} = -(x+\sigma \eta_t)^3
% \end{align*}
% is equivalent to the It\^o SDE 
for which the direct steady-state integral calculation indeed diverges.
Enforcing relation \ref{eq:consistency_Ito} to leading order in $\D{t}$ for this system gives 
\begin{align}
%    0 &= \int\limits_{-\infty}^{\infty}\rho^*(x)  \left[2xF(x) + G(x)^2\right] \dd{x} \nonumber\\
    0 = &\int\limits_{-\infty}^{\infty}\rho^*(x) 
        \big[2x(-x^3-3\sigma^2 x) \nonumber \\
        & \qquad \qquad \quad + (9\sigma^2x^4 + 36 \sigma^4 x^2 + 15 \sigma^6) \big] \dd{x}\nonumber\\
    = &15 \sigma^6\int\limits_{-\infty}^{\infty}\rho^*(x) \dd{x} + (36\sigma^4-6\sigma^2)\int\limits_{-\infty}^{\infty}x^2\rho^*(x) \dd{x} \nonumber\\
    & \qquad \qquad \quad + (9\sigma^2-2)\int\limits_{-\infty}^{\infty}x^4\rho^*(x)\dd{x}\;. \nonumber\\
%    = &15 \sigma^6  +(36\sigma^4-6\sigma^2)\mu_2 + (9\sigma^2-2)\mu_4 \nonumber \\
    = &15 \sigma^6  +6\sigma^2(6\sigma^2-1)\mu_2 + (9\sigma^2-2)\mu_4 
    \label{eq:surface}
 \end{align}
So we obtain a relationship between moments of the equilibrium $\rho^*$. 

However we notice a problem: if $\sigma$ is large enough that $9\sigma^2-2>0$ and $6\sigma^2-1>0$ (i.e., $\sigma > \sqrt{2}/3$), all terms on the right hand side are positive and there is no way for the equality to hold.

If we had preserved all terms from Eq.~\eqref{eq:consistency_Ito}, rather than truncating at leading order, we would have obtained the full, exact relation
\begin{align} 
    0 &=  15 \sigma^6  +(36\sigma^4-6\sigma^2+9\sigma^4 \D{t})\mu_2 \nonumber \\
    & \qquad\qquad + (9\sigma^2-2 +6 \sigma^2 \D{t})\mu_4 +  \D{t}\  \mu_6\;.
\end{align}
This still does not avoid the problematic implication at large $\sigma$---in fact, it makes the situation slightly ``worse'' by adding more positive terms.  This contradiction implies that we must have been wrong to treat $\mu_2$ as finite (implicit in utilizing relation \ref{eq:consistency_Ito})---i.e., the equilibria for these values of $\sigma$ must have divergent second moments. %\footnote{We note that the earlier assumption of finite underlying stochastic-process variance was only necessary to arrive at the It\^o SDE from the initial Langevin-type equation; our analysis of the It\^o SDE itself does not rely on that assumption.}

\subsection{Higher-Moment Methods}
If we repeat our above analysis, but with the $2k^\textrm{th}$ raw moment of $\rho^*$ instead of the second\footnote{For symmetric equilibria like our cubic example, odd moments are all zero.}, we have
\begin{align}
&\mu_{2k} = \int\limits_{-\infty}^{\infty}\xi^{2k}\rho^*(\xi) \dd{\xi} \nonumber \\
    & = \int\limits_{-\infty}^{\infty}\rho^*(x) 
       \int\limits_{-\infty}^{\infty}
         \frac{\xi^{2k}}{G(x)\sqrt{2 \pi \D{t}}}e^{\frac{-[\xi - x - F(x)\D{t}]^2}{2 G(x)^2 \D{t}}} \dd{\xi}\ \dd{x}\;. %\nonumber \\
    % & = \int\limits_{-\infty}^{\infty}\rho^*(x) \left\{ I_{2k}[x+F(x)\D{t},G(x)\sqrt{\D{t}} ] \right\}     
\end{align}
Integrals of the following form arise:
\begin{align*}
    I_{2k} :&= \frac{1}{\sigma \sqrt{2 \pi}} \int\limits_{-\infty}^{\infty}u^{2k} e^{\frac{-(u-a)^2}{2\sigma^2}} \dd{u}\\
%     &= \frac{1}{\sigma \sqrt{2 \pi}} \int\limits_{-\infty}^{\infty}(w+a)^{2k} e^{\frac{-w^2}{2\sigma^2}} \dd{w} \\
%     &= \sum_{i = 0}^{k} {2k \choose 2i} a^{2k-2i} \frac{1}{\sigma \sqrt{2 \pi}} \int\limits_{-\infty}^{\infty}w^{2i} e^{\frac{-w^2}{2\sigma^2}} \dd{w}
% \end{align*}
% \dma{to do: remove the work here and just state the definition of I2k and the result.}
% since all the odd terms drop out. This gives
% \begin{align*}
%     I_{2k} &= \sum_{i = 0}^{k} {2k \choose 2i} a^{2k-2i} (2i-1)!!\sigma^{2i} \\
    &= (2k)! \sum_{i = 0}^{k} \frac{\sigma^{2i}a^{2k-2i}}{(2i)!!(2k-2i)!} \;.
\end{align*}
%\dma{Not a big deal but it would be nicer to use $l$ instead of $i$ to avoid confusion with imaginary numbers.}
So with any It\^o SDE we have
\begin{align*}
    &\mu_{2k} = \int\limits_{-\infty}^{\infty}\rho^*(x) \\
    & \times \left[ (2k)! \sum_{i = 0}^{k} \frac{(G(x)\sqrt{\D{t}})^{2i}(x+F(x)\D{t})^{2k-2i}}{(2i)!!(2k-2i)!} \right] \dd{x} \\
    &= \sum_{i=0}^{k}  \frac{ (2k)!}{(2i)!!(2k-2i)!} \D{t}^i\int\limits_{-\infty}^{\infty}\rho^*(x) G(x)^{2i} \\
    & \qquad\quad\quad \times \sum_{j=0}^{2k-2i} {2k-2i \choose j} x^j [F(x)\D{t}]^{2k-2i-j} \dd{x} \;.
\end{align*}
Regrouping by powers of $\D{t}$ and retaining only leading order behavior, we find that the constant term ($i=0, j=2k$) cancels from the left hand side, leaving
\begin{align}
    0 &= 
    \D{t}\int\limits_{-\infty}^{\infty}\rho^*(x)\left[2kx^{2k-1}F(x) + G(x)^2  \right] \dd{x} \;.
%    0 &= 
%    \D{t}\int\limits_{-\infty}^{\infty}\rho^*(x)\left[2kx^{2k-1}F(x) + G(x)^2  \right] \dd{x} + O(\D{t}^2) \;.
    % & \ \ + \D{t}^2\Bigg\{\int\limits_{-\infty}^{\infty}\rho^*(x)\big[ x^{2k-2}F(x)^2 \\
    % & \qquad \quad + G(x)^2 x^{2k-1}F(x) + G(x)^4  \big]\dd{x}\Bigg\} + O(\D{t}^3) \\
    % &\implies \int\limits_{-\infty}^{\infty}\rho^*(x)\left[x^{2k-1}F(x) + G(x)^2 \dd{x} \right] \sim O(\D{t}) \\
    % & \quad \sim \D{t}\left\{\int\limits_{-\infty}^{\infty}\rho^*(x)\left[ x^{2k-2}F(x)^2 + G(x)^2 x^{2k-1}F(x) + G(x)^4  \right]\dd{x}\right\} + O(\D{t}^2)
    \label{eq:General_result}
\end{align}
%const: i=0, j=2k
%dt: (i=0, j=2k-1) and (i=1, j=2k)
%dt^2: (i=0, j=2k-2) and (i=1, j=2k-1) and (i=2, j=2k)

% For this (leading order) equality to hold, we should have (for $k=1,2,\ldots$)
% \begin{align}
%     0&= 
%   \int\limits_{-\infty}^{\infty}\rho^*(x)\left[2kx^{2k-1}F(x) + G(x)^2 \right] \dd{x} .
% \end{align}

This relation should hold for any equilibrium of an It\^o SDE for which the $2k^{\textrm{th}}$ raw moment is finite. If $F(x)$ and $G(x)^2$ are polynomials, this may be used to obtain a recursion relation for all moments of the equilibrium $\rho^*$.

For example, in the case of the cubic generalized-Langevin attractor from \cite{sabinmiller2024interpretation},
\begin{align}
    \frac{\dd{x}}{\dd{t}} &= -X^3, \textrm{ where } X \sim N(x,\sigma) \nonumber\\
    \implies \dd{x} &= (-x^3-3\sigma^2 x)\ \dd{t}  \nonumber \\
    & \qquad + \sqrt{15\sigma^6 + 36 \sigma^4 x^2 +9 \sigma^2 x^4}\ \dd{W}, \label{eq:nIto_form_of_cubic}
\end{align}
(equivalence logic argued in that paper: $F$ = mean and $G$ = standard deviation of the d$x$/d$t$ distribution) would leave us with the equilibrium relation
\begin{align*}
    0 = &\int\limits_{-\infty}^{\infty}\rho^*(x)\Big[2kx^{2k-1}(-x^3-3\sigma^2x) \\
   &\qquad + (9\sigma^2x^4+ 36\sigma^4x^2+ 15\sigma^6) \Big] \dd{x} \\
   = &15 \sigma^6 + 36 \sigma^4 \mu_2 + 9 \sigma^2\mu_4 -6k\sigma^2 \mu_{2k} -2k\mu_{2k+2}
\end{align*}
for integers $k\geq 1$. %Unfortunately this gives a slightly under-specified system of equations for the moments $\mu_{2k}$ in terms of the inherent noise amplitude $\sigma$.
% :
% \begin{equation*}
%     % 0 &= 15 \sigma^6 + (36\sigma^4-6\sigma^2)\mu_2 + (9\sigma^2-2)\mu_4\\
%     % 0 &= 15 \sigma^6 + 36\sigma^4\mu_2 + (9\sigma^2-12\sigma^2)\mu_4 -4\mu_6 \\
%     0 = 15 \sigma^6 + 36\sigma^4\mu_2 + 9\sigma^2\mu_4 -6k\sigma^2\mu_{2k}-2k\mu_{2k+2} \;, %, \\
%     %& \qquad\qquad\qquad\qquad\qquad\qquad\qquad\qquad\qquad\quad k\geq1
%     % \label{eq:moment_recursion}
% \end{equation*}
%valid for integers $k \geq 1$.

\begin{figure}[th!] 
\centering
	\includegraphics[width=\columnwidth]{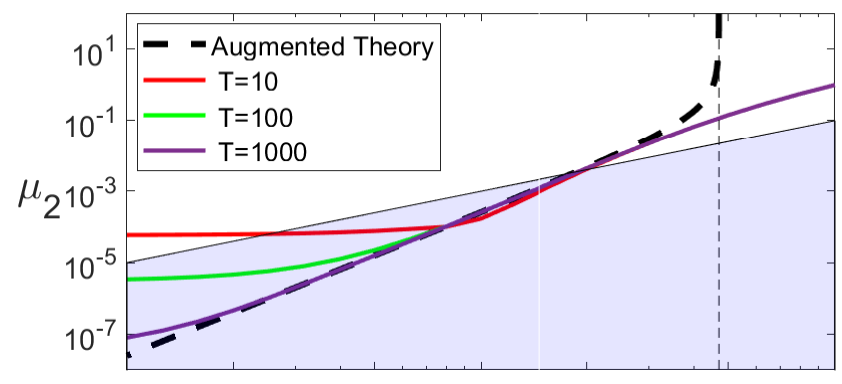}
	\includegraphics[width=\columnwidth]{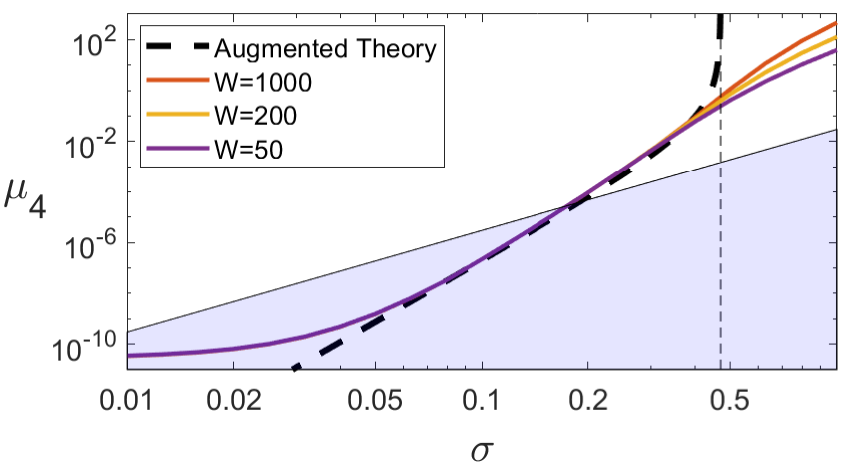}
	\caption{\textbf{Numerical validation.}  Comparison of numerical results (via Fokker-Planck evolution) to the theoretical relation, augmented with the extra Gaussian condition $\mu_4 = 3 \mu_2^2$. The shaded region indicates $\mu_2 < 0.1 \sigma^2$ (top) and correspondingly $\mu_4 < 0.03 \sigma^4$ (bottom), where the Gaussian approximation (from $\mu_2 \ll \sigma^2$) should be most valid. \textbf{Top:}   Smaller $\sigma$ values take longer simulated time $T$ to equilibrate, but do approach the theorized line in the shaded region.  For high noise amplitudes, the relation need not hold, and indeed theory suggests that $\mu_2$ and $\mu_4$ should diverge for $\sigma>\sqrt{2}/3$ (indicated by the vertical dashed line)---though this divergence is invisible at constant domain width. \textbf{Bottom:} As predicted by theory, the fourth moment $\mu_4$ does indeed appear to diverge for $\sigma>\sqrt{2}/3$, though simulation with ever wider domain width $W$ (measured in number of standard deviations of the equilibrium solution) is needed capture more of the distribution's tails (all curves shown for $T = 100$).
	%\dma{Plots not really readable for colorblind, but we can make alternative version during review. } 
	} 
	\label{fig:2d_matching}
\end{figure}

While this slightly under-specified system of equations doesn't yield exact moments, it implies that those moments should lie on a surface, which we confirm by numerical simulation (see section S2 of the SM).  When the typical magnitude of $x$ is small compared to $\sigma$ (i.e., $\mu_2 \ll \sigma^2$), however, Eq.~\eqref{eq:nIto_form_of_cubic} is well approximated by an SDE with constant noise and linear drift: an Ornstein–Uhlenbeck process \cite{uhlenbeck1930theory, wang1945theory} (see also, e.g., \cite{van1992stochastic} or \cite{gillespie1991markov}). This implies a normal distribution at equilibrium, with moment relationship 
\begin{equation}
    \mu_4 = 3 \mu_2^{2}.
\end{equation}
Plugging this additional constraint into our lowest-order relation Eq.~\eqref{eq:surface} yields 
\begin{equation*}
    0 = 15 \sigma^6 + 6 \sigma^2(6\sigma^2 - 1)\mu_2 + 3 (9\sigma^2-2)\mu_2^2  \;, %, 
\end{equation*}
which agrees well with simulation in the relevant parameter region: see Fig.~\ref{fig:2d_matching}. For a direct look at the Gaussian nature of equilibria across this transition, see section S3 of the SM.

% \begin{figure}[th!] 
% \centering
% 	\includegraphics[width=\columnwidth]{matching_theory_surface.eps}
% 	\caption{\textbf{Numerical validation.}  Fokker-Planck solutions to the ``nIt\^o''-converted cubic system with different values of $\sigma$ confirm that the equilibrium does obey the leading-order moment-relation equation \eqref{eq:surface}. 
% \dsm{do orthogonal distance, lagrange multipliers, but then no relative error... unless (distance to surface) / (distance to origin)}
% \dsm{wider aspect ratio perspective?}
% \dsm{show divergence with larger sigma?}}
% 	\label{fig:surface}
% \end{figure}

%\dma{Maybe a single multi-panel figure showing the correctness of moment relationships, and also showing the divergence at critical sigma.  Or maybe multiple panels aren't necessary?}
\section{Conclusions}
% The first proposition of this paper---the argument for Ito-equivalency of nonlinear Langevin-type systems---is really a proposed definition rather than a theoretical result.  Like Langevin equations themselves, the notation is simple and intuitive, but solid mathematical interpretation requires the use of the more rigorous notation, and we propose that interpretation in terms of It\^o calculus. 

% We apply logic based on the central limit theorem for finite-variance random variables, but the Langevin noise terms are not regular random variables and their variance may not be well-defined or finite. If variance is treated as well-defined but not finite, other (non-Gaussian) stable distributions per time-step may arise, rather than normally distributed It\^o time-steps. 
% %\dma{revise a bit: not exploratory but an interpretation due to the limits of Langevin's notation.}

% We also note the perhaps-undesirable sensitivity to the assumption of Gaussian underlying noise in Eq.~\eqref{eq:generalsys}.  
% %\sout{Due to integration over the ``Langevin-style'' noise, the resulting It\^o system will depend on the shape of that underlying distribution.} 
% In particular, the assumption that $\eta_t$ is normally distributed may be incorrect for some systems with biased or irregularly shaped noise, and if the noise shape is known it should be used. %However, a multiple-time-scale argument may allow the Central Limit Theorem to apply twice

%Appropriateness of CLT for non-standard ``random variable", and CLT Gaussian -> Ito?
We have introduced an equilibrium-analysis technique for It\^o SDEs inspired by the observation that such equilibria should theoretically be fixed under sufficiently precise numerical integration. The resulting analysis, which culminates in Eq.~\eqref{eq:consistency_Ito}, should apply to any It\^o system with an equilibrium where the second raw moment of that equilibrium is finite, but it is of particular use when the functions $F$ and $G^2$ are polynomial in nature, since this allows the analysis to yield explicit relations between even moments of the equilibrium rather than merely integrals against that unknown distribution. 

We applied this technique to an example system arising from previous work, and were able to prove (by contradiction) that this system has a critical value of the noise parameter $\sigma$ at which its equilibrium must have divergent moments, since the relation becomes contradictory. This result would be very difficult to glean from direct numerical integration of the It\^o (or corresponding Fokker-Planck) equation in question, due to the subtlety of this divergence in the distribution's tails for any finite domain width.

In the analysis of this example system in the parameter region where its equilibrium has finite moments, we utilized an additional, near-Gaussian approximation (valid when $\mu_2 \ll \sigma^2$) which allowed us to fully prescribe the moments of the equilibrium as a function of the noise parameter $\sigma$, and confirmed with numerical integration that the relation appears to hold---simulations converged to this augmented relation in the region of this approximation's validity in forwards time. It remains unclear whether a more general constraint valid for arbitrary $\sigma$ closer to the moment-divergence boundary can be found. If such an additional constraint could be found and to the extent that it generalizes, this technique could enable the full specification of equilibrium moments. Regardless, we hope this analysis may yield some useful insight into specific systems which resist other techniques.

\begin{center}
  \pgfornament[width = 0.8\columnwidth, color = black]{88}
\end{center}

\begin{acknowledgments}
The authors thank Bill Kath for useful conversation, Gary Nave for help with relevant literature, and the National Science Foundation for support through the Graduate Research Fellowship Program.
\end{acknowledgments}

% \bibliography{refs}


\pagebreak

\pagebreak 
\newpage 
\clearpage



\section*{Supplementary material}

# Supplemental Material: Equilibrium Moment Analysis of Itô SDEs

David Sabin-Miller*

*Center for the Study of Complex Systems, University of Michigan, Ann Arbor, MI, USA*

Daniel M. Abrams

*Department of Engineering Sciences and Applied Mathematics,
Northwestern University, Evanston, IL, USA
Northwestern Institute for Complex Systems,
Northwestern University, Evanston, IL, USA and
Department of Physics and Astronomy, Northwestern University, Evanston, IL, USA*

## S1. ADDITIONAL FIGURES

Here we have provided some additional 3-dimensional figures elaborating on the 2-dimensional ones from the main text.

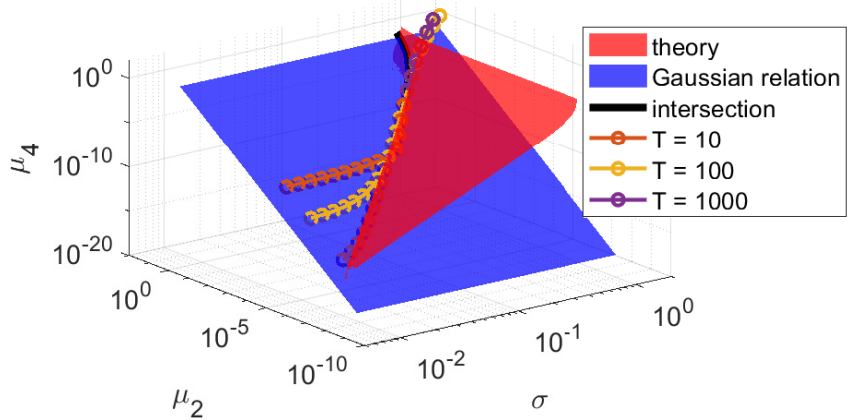

FIG. S1: **Small-noise convergence.** Three dimensional version of Fig. 1 (top panel) from the main paper; at small noise $\sigma$, solutions converge over increasing simulated time $T$ to the intersection of the Gaussian condition and Eq. (11), our theorized surface from the main paper relating the equilibrium's second and fourth moments to the inherent noise $\sigma$: $0 = 15\sigma^6 + 6\sigma^2(6\sigma^2 - 1)\mu_2 + (9\sigma^2 - 2)\mu_4$.

---

* dasami@umich.edu

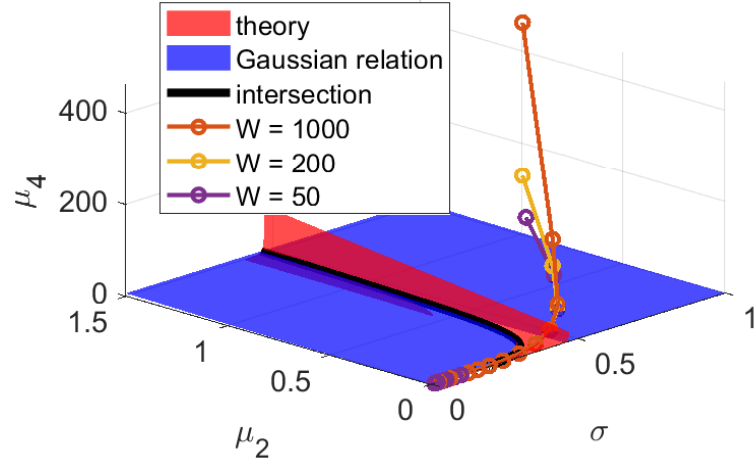

FIG. S2: **Large-noise divergence.** Three dimensional version of Fig. 1 (bottom panel) from main paper relating noise amplitude $\sigma$ to the second and fourth moments of the equilibrium ($\mu_2$ and $\mu_4$, respectively), but on a linear scale. We can see the clear suggestion of divergence for (at least) $\mu_4$; as our domain captures more standard deviations $W$ of the solution, the measured $\mu_4$ grows without bound. We propose that noise values $\sigma$ beyond the asymptote "curtain" $\sigma^* = \sqrt{2}/3 \approx 0.47$ must have divergent second and fourth moments.

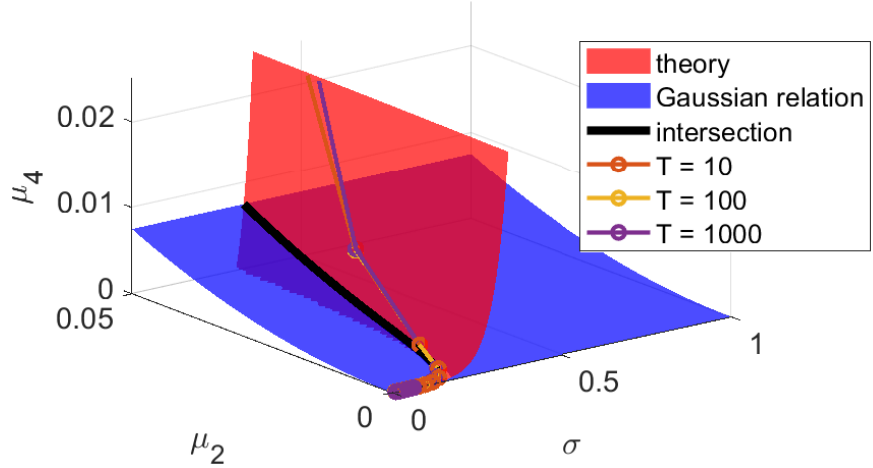

FIG. S3: **Gaussian validity boundary.** Three dimensional, linear-scale zoom of the region where the Gaussian assumption fails to hold. Our solutions fall off the intersection line due to their non-Gaussian nature.

## S2. GAUSSIAN/NON-GAUSSIAN TRANSITION

In the main text we argue that the small-noise limit of the stochastic-cubic-attractor system leads to a Gaussian equilibrium. We can see this transition clearly in Fig. S4.

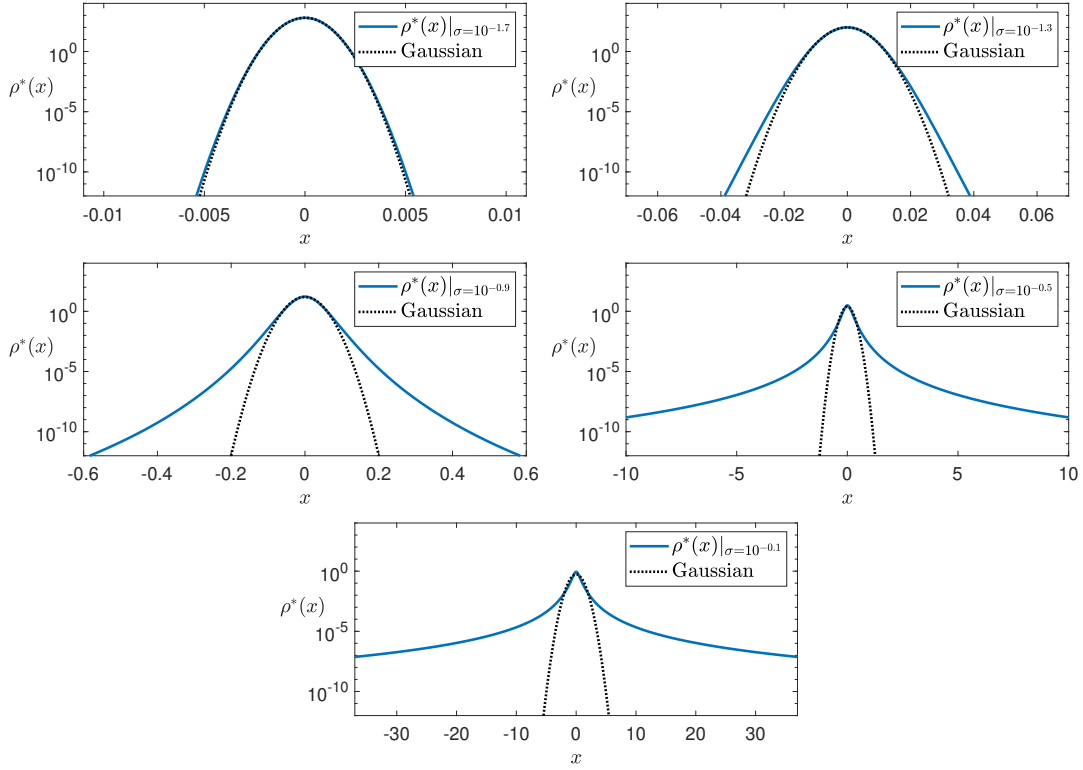

FIG. S4: **Equilibrium shape transition.** Comparison of equilibrium distributions for different noise values to Gaussian distributions with the same standard deviation. **Top Left:** For small $\sigma$, equilibria exhibit the signature parabolic shape (on semi-log axes) indicating a near-Gaussian distribution. **Other panels:** For larger $\sigma$, equilibria exhibit increasingly "fat tails" differentiating them from Gaussians.


% \newpage 
% \clearpage
% \includepdf[pages={4}]{Sabin-Miller_SM.pdf}
% \newpage 
% \clearpage
% \includepdf[pages={5}]{Sabin-Miller_SM.pdf}
% \newpage 
% \clearpage
% \includepdf[pages={6}]{Sabin-Miller_SM.pdf}

\end{document}